\documentclass{gtpart}     
\usepackage{pinlabel}  
\usepackage{amscd}
\author[Patrick Gilmer]{ Patrick M. Gilmer}
\address{Department of Mathematics\\
Louisiana State University\\
Baton Rouge, LA 70803\\
USA}
\email{gilmer@math.lsu.edu}
\urladdr{http://www.math.lsu.edu/\textasciitilde gilmer/}

  \title[Congruence and Quantum Invariants ]{Congruence and Quantum Invariants of 3-manifolds}

\arxivreference{math.GT/0612282}

\volumenumber{}
\issuenumber{}
\publicationyear{}
\papernumber{}
\startpage{}
\endpage{}
\doi{}
\MR{}
\Zbl{}
\received{}
\revised{}
\accepted{}
\published{}
\publishedonline{}
\proposed{}
\seconded{}
\corresponding{}
\editor{}
\version{}

\newtheorem{thm}{Theorem}[section]
\newtheorem{lem}[thm]{Lemma}
\newtheorem{prop}[thm]{Proposition}
\newtheorem{cor}[thm]{Corollary} 
\newtheorem{de}[thm]{Definition} 
\newtheorem{rem}[thm]{Remark}

\newcommand{\BH}{{\mathbb{H}}}
\newcommand{\BZ}{{\mathbb{Z}}}
\newcommand{\BQ}{{\mathbb{Q}}}

\newcommand{\BO}{{\mathcal{O}}}
\newcommand{\BD}{{\mathbb{D}}}

\newcommand{\Si}{{\Sigma}}

\newcommand{\bb}{{\mathfrak{b}}}

\newcommand{\BC}{{\mathbb{C}}}

\DeclareMathOperator{\Int}{Int}

\begin{document} 

\begin{abstract}
 Let $f$ be an integer greater than one. We study three progressively finer equivalence relations on
closed 3-manifolds generated by Dehn surgery with denominator $f$: weak $f$-congruence,  $f$-congruence, and strong $f$-congruence. If $f$ is odd, weak $f$-congruence preserves the ring structure on cohomology with $\BZ_f$-coefficients.  We show that strong  $f$-congruence coincides with a relation previously studied by Lackenby. 
Lackenby showed that the quantum $SU(2)$ invariants are well-behaved under this congruence. We strengthen this result and extend it to the $SO(3)$ quantum invariants. We also obtain some corresponding results for the coarser equivalence relations, and for quantum invariants associated to more general modular categories. We compare $S^3$, the Poincar\'{e} homology sphere, the Brieskorn homology sphere
 $ \Si(2,3,7)$ and their mirror images up to  strong  $f$-congruence. We distinguish the weak $f$-congruence classes of some manifolds with the same $\BZ_f$-cohomology ring structure.  
 \end{abstract}

\maketitle


 \section{Introduction}
 
 
 Weak type-$f$ surgery is a kind of surgery along a knot in a 3-manifold which generalizes the notion of $n/sf$ surgery in a homology sphere. Such surgeries preserves the cohomology groups with $\BZ_f$-coefficients. Weak type-$f$ surgery generates an equivalence relation on  3-manifolds which we call  weak 
 $f$-congruence. If $f$ is odd, we show that weak type-$f$ surgery also preserves  the cohomology ring structure with $\BZ_f$-coefficients.  
 
Strong type-$f$ surgery is a kind of surgery along a knot in a 3-manifold which generalizes the notion of $1/sf$ surgery in a homology sphere.
We call the equivalence relation on 3-manifolds that it generates strong $f$-congruence. Motivated by Fox's notion of $f$-congruence for links \cite{F}, Lackenby defined an equivalence relation on 3-manifolds which  he called congruence modulo $f$ \cite{La}.   Congruence modulo $f$ is generated  by a move which increments the framing on a component of framed link description  by $f$.  We show that strong $f$-congruence coincides with congruence modulo $f$.
We also consider a relation which we call $f$-congruence. It is coarser than strong $f$-congruence and finer than  weak $f$-congruence.  
 
 We refine a relation, due to  Lackenby,  between the $SU(2)$ quantum invariants of  manifolds which are congruent modulo $f$. We find relations reflecting congruence and weak  congruence. We also discuss quantum invariants associated to some modular categories. 
 
{\it  In this paper,  $p$ will denote  an odd prime.} We show that the quantum $SO(3)$ invariant at a $p$th root of unity  is preserved, up to phase, by $p$-congruence.  We  give  finite lists of the only possible $f$ for which there might be strong $f$-congruences between $S^3$, the Brieskorn homology spheres $\pm \Si(2,3,5)$ and $\pm \Si(2,3,7)$. Moreover we realize some of these strong congruences.  
   
   We also show that the quantum $SO(3)$ invariant at a $p$th root of unity has a very simple surgery formula for  weak type-$p$ surgeries. As a corollary, we distinguish, up to weak $p$-congruence  for all $p>3$,  two  3-manifolds with the  same  cohomology rings: 0-framed surgery to the Whitehead link and $\#^2 S^1 \times S^2$.  This can also be done using work of Dabkowski and  Przytycki's \cite[Theorem 2(i)]{DP} on Burnside groups of links. We also distinguish another manifold from $\#^2 S^1 \times S^2$ up to weak $p$-congruence for small $p.$

 We strengthen a result of Masbaum and the  author on the divisibility of certain quantum invariants.

We thank Gregor Masbaum, Brendan Owens, Jozef Przytycki, Khaled Qazaqzeh and  several referees for comments, suggestions and/or discussions. This research was partially supported by NSF-DMS-0604580.

 \section{Congruence}
 
 Our convention is that all manifolds are compact, and oriented, unless they fail to be compact by construction. We use a minus sign to indicate orientation reversal.
We use $N$, $N'$  and $M$ to  denote closed connected 3-manifolds.
 {\it  In this paper, we  let $f$ denote an integer greater than one.}

\begin{de}[Lackenby] Two closed 3-manifolds are congruent modulo $f$ if and only if they possess framed link diagrams which are related by a sequence of moves:  the usual Kirby moves and also the move  of changing the framings by adding  multiples of  $f$. \end{de}

Suppose $\gamma$ is a simple closed curve in a closed connected 3-manifold $N$.  Let $\nu_\gamma$ denote a closed tubular neighborhood of $\gamma$, and $ \mathcal {T}_\gamma$ the boundary of $\nu_\gamma$. By a meridian for $\gamma$, we mean a simple closed curve $\mu$ in $ {\mathcal T}_\gamma$ which is the boundary of a transverse disk to $\gamma$. By a longitude, we mean a simple closed curve $\lambda$ in $ {\mathcal T}_\gamma$ which meets a meridian in a single point transversely.
The process of removing $\nu_\gamma$ from $N$ and reattaching it so that a curve $\mu'$  that is homologous to  $ \ell  \lambda+ n \mu$ bounds a disk in the reglued solid torus will be called an $n/\ell$
surgery to 3-manifold $N$ along $\gamma$. Here $n$, $\ell$ are integers, and $n$ is relatively prime to $\ell.$ The denominator of the surgery is $\ell$. This is well defined (i.e. does not depend on the choice of $\lambda$) up to sign. The congruence class of $n$ modulo $\ell$  is well defined up to sign, and $n$ is called the numerator for the surgery.

\begin{de}
A $n/\ell$ surgery is called weak type-$f$ surgery if $\ell \equiv 0 \pmod{f}$.  A  $n/\ell$-surgery is called type-$f$ surgery if  $\ell \equiv 0 \pmod{f}$ and  $n \equiv \pm1 \pmod{f}$.  
A  $n/\ell$-surgery is called strong type-$f$ surgery if  $\ell \equiv 0 \pmod{f}$ and  $n \equiv \pm1 \pmod{\ell}$. 
\end{de}

 If we may obtain $N'$ from $N$ by a strong, weak or plain type-$f$ surgery, we may also obtain $N$ from $N'$ by a  type-$f$ surgery of the same variety (simply by reversing the process).

\begin{de} The equivalence relation generated by  strong type-$f$ surgeries is called strong  
$f$-congruence. 
The equivalence relation generated by   type-$f$ surgeries is called $f$-congruence.  The even coarser equivalence relation generated by weak type-$f$ surgeries will be called weak $f$-congruence.\end{de}

 \begin{prop}\label{obv} Let $m$ be a positive integer.
If $M$ is (respectively weakly, strongly) $fm$-congruent  to $N$,
then $M$ is (respectively weakly, strongly) $f$-congruent to $N$. 
\end{prop}

\begin{thm}\label{ce} Two 3-manifolds are congruent modulo $f$ if and only if they are  strongly $f$-congruent . \end{thm}

\begin{figure}[ht!] 
   \labellist
    \small
  \hair 2pt
 \pinlabel {$ \begin{array} {c} \text{insert a surgery} \\ \text{curve with a}\\ \text{$0$-framed meridian} \end{array}$ }  at  1 70
 \pinlabel {$n$} at  90.3 108.5
 \pinlabel {$0$} at  115 57
  \pinlabel {$ \begin{array} {c} \text{congruence}\\ \text{modulo $f$} \\ \text{move} \end{array}$ } at  165 70
  \pinlabel {$n$} at   230 108.5
   \pinlabel {$-sf$} at  260 57
     \pinlabel {$\text{slam dunk}$}  at  285 80
       \pinlabel { $\frac {1+sfn}{sf}$} at  365 66
  \endlabellist
  \begin{center}
  \includegraphics[width=3in]{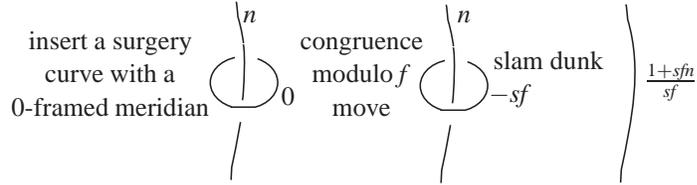}
 \caption{ strong $f$-congruence move generates 1/f surgery 
( $n- \frac 1 {-sf}= \frac {1+fn}{sf}$)}  \label{f1}
\end{center}
\end{figure}

\begin{figure}[ht!] 
 \labellist
    \small
  \hair 2pt
 \pinlabel {$n$} at  11 115
   \pinlabel {$ \begin{array} {c} \text{do}\\  {-1/f} \\ \text{surgery} \end{array}$ } at  42 60
  \pinlabel {$n$} at   120 115
   \pinlabel {$-1/f$} at  160 55
     \pinlabel {$\text{Rolfsen twist}$} at  215 80
       \pinlabel { ${n+f}$} at  310 110  \endlabellist
\begin{center}
\includegraphics[width=2.5in]{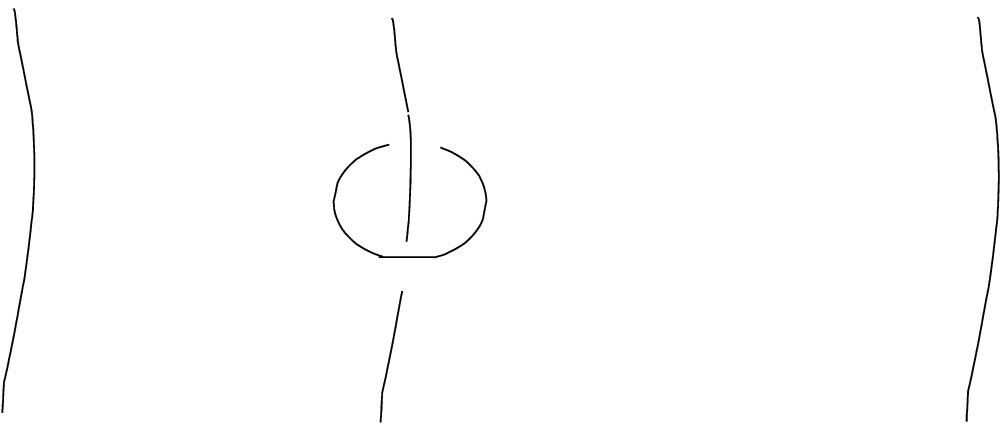}
\caption{1/f surgery  generates strong $f$-congruence move}\label{f2}

\end{center}
\end{figure}

 \begin{proof} Suppose $N$ is already described by surgery on a link $L$ in $S^3$. We want to see that the result of  $\frac{1+sfn}{sf}$ (for any $s \in \BZ$) to $N$ along a knot $K$ in the complement of $L$ is strongly $f$-congruent to $N$. By a well-known trick \cite[Prop 5.1.4]{GS}, we can get another surgery  description of $N$ by including $K$ with framing $n$ and a meridian of $K$ framed zero. 
We can then  change the framing on the meridian from zero to $-sf$  \'{a} la  Lackenby, then we may do a slam dunk \cite[p.163]{GS}. See Figure \ref{f1}. The result of  $\frac{-1+sfn}{sf}$  to $N$ can be realized similarly.

Suppose now $N$ is already described by surgery on a framed link $S^3$ which includes a  component $K$ framed, say $n$. If we perform $-1/f$ surgery on a meridian of $K$, and then do a Rolfsen twist \cite[p.162]{GS}, we will  have changed the framing on $K$ to $n+f$. See Figure \ref{f2}.
 \end{proof} 
 
We  need the  concept of an $f$-surface for the next proof. This concept  is also required to formulate some  later results. The idea here is that of a generalized surface where a number of sheets which is multiple of $f$ are allowed to coalesce along  circles.  Note that  a non-orientable closed surface together with a selected one manifold dual to the first Stiefel-Whitney class of  the surface  and a choice of orientation on the complement of this one manifold is a simple example  of a good $2$-surface.

\begin{de} An $f$-surface $F$ is the result of attaching, by a map $q$, the whole
boundary of an oriented surface $\hat F$ to a collection of circles $\{ S_i\}$ by a map  which when restricted to the  inverse image under $q$ of each $S_i$ is a $f  t_i $-fold ( possibly disconnected)  covering space of $S_i$. If each component of each $q^{-1} S_i$ is itself a covering space of $S_i$ with degree divisible by $f$, we say $F$ is a good $f$-surface. The image of the interior of the surface is called the 2-strata. The image of the boundary is called the 1-strata. If only part of the boundary of $F$ is so attached, we call this a $f$-surface with boundary, and the image of the unattached boundary is called the boundary. \end{de}

\begin{thm}\label{ring} A weak $f$-congruence between $N$ and $N'$ induces a graded group isomorphism between  $H_*(N,\BZ_f)$ and $H_*(N',\BZ_f)$ and between $H^*(N,\BZ_f)$ and $H^*(N',\BZ_f)$. If $f$ is odd, this induced  isomorphism preserves the ring structure.  If $f$ is two, this induced  isomorphism need not preserves the ring structure.
\end{thm}

\begin{proof} Let $\gamma$ denote the curve in $N$ that we perform the weak type-$f$ surgery along.  Let $ \gamma'$ denote the core of the new solid tori in $N'$. Let $X=N\setminus \gamma =  N'\setminus \gamma'$. 
As the maps $H_1(X,\BZ_f) \rightarrow H_1(N,\BZ_f)$ and  $H_1(X,\BZ_f) \rightarrow H_1(N', \BZ_f)$ induced by the inclusions are surjective and have the same kernel, it follows that  the induced mappings$:H^1(N, \BZ_f) \rightarrow H^1(X,\BZ_f)$  and $:H^1(N',\BZ_f) \rightarrow H^1(X,\BZ_f)$ are injective and have the same image. Thus these maps induce isomorphisms on $H^1( \ , \BZ_p)$
and $H_1( \ , \BZ_p)$. Connectivity yields isomorphisms on $H_0( \ , \BZ_p)$, and orientations  yields isomorphisms on $H_3( \ , \BZ_p).$
The above isomorphisms and Poincar\'{e} duality  yield the others.

Using Poincar\'{e} duality and the equation $(a \cup b )\cap z= a  \cap (b  \cap z)$, to see that the ring structure is preserved it suffices to check that the isomorphism on $H^1( \ , \BZ_f)$  preserves the trilinear triple product $(\chi_{1} \cup \chi_{2} \cup \chi_{3})\cap [N].$

To verify this, we use $f$-surfaces 
 to represent classes in 
$H^1(\ ,\BZ_f)$.
An $f$-surface has a fundamental class $H_2(F, \BZ_f)$ which is given by  the sum of the oriented 2-simplices in a triangulation of $\hat F$. Thus a  $f$-surface $F$ embedded in $N$ represents an element $[F] \in H_2(N, \BZ_f).$ Poincar\'{e} dual to $[F]$ is the cohomology class $\chi_F \in H^1(N, \BZ_f)$  which may also be described by the (signed)
intersection number of a loop which meets  $F$ transversely in the 
2- strata. Every 
cohomology class in $H^1(N, \BZ_f)$ may be realized in this way by an $f$-surface.

Given an $f$-surface $F$ in $N$, we may isotope $F$ so  it  transversely intersects $\gamma$ in the 2-strata.   Each circle component of $F \cap \mathcal{T}_\gamma$ consists of a collection  of meridians of $\nu_\gamma$. Viewing $F \cap \mathcal{T}_\gamma$ from the point of view of ${\gamma'}$ we see a collection of parallel torus knots.  A component  is homologous   to a number of longitudes  of  $\gamma'$ which is divisible by $f$ plus some number of meridians of  $\gamma'$ (necessarily prime to $f$).  
Thus $F \setminus (F \cap \nu_\gamma)$ maybe completed to an $f$-surface by adjoining the mapping cylinder of the projection of $(F \cap \mathcal{T}_\gamma)$ to  $\gamma'.$  Let $F'$ denote the new $f$-surfaces in $N'$  constructed in this manner. The induced  isomorphism from $H^1(N, \BZ_f)$ to $H^1(N', \BZ_f)$  sends $\chi_F$ to $\chi_F'$.
We could also complete $F \setminus (F \cap \nu_\gamma)$ to  form $F'$ in some other way by adding any $f$ surface with boundary in $\nu_\gamma'$ with boundary $(F \cap \mathcal{T}_\gamma)$. The class of $[F']$ does not depend on this choice, as $H_3(\nu_\gamma',\BZ_f)$ is zero.

Any three $f$-surfaces  $F_1$, $F_2$ and $F_3$ in $N$  may be isotoped  so that $F_1 \cap F_2 \cap F_3$ lies in the intersection of the 2-strata of  these surfaces and consists of a finite number points and, in a neighborhood of these triple points,  the three surfaces  look locally like the intersection of the three coordinate planes in 3-space. One has that the triple product $(\chi_{F_1 } \cup \chi_{F_2} \cup \chi_{F_3})\cap [N]$ can be computed as the number of triple points as above counted according to sign in the usual manner, and denoted
  $F_1 \cdot F_2 \cdot F_3$.  This number only depends on the homology classes:
  $[F_1],$$[F_2]$,$[F_3].$ 
  Note that  $F_1 \cdot F_3 \cdot F_2
  = - F_1 \cdot F_2 \cdot F_3.$

Given $F_1$, ${F_2}$ and  ${F_3}$ in $N$, we can isotope them so that the  intersections  of $\gamma$ with the $F_i$  are all grouped together as one travels along $\gamma$, and that the intersections are encountered first with (say) $F_1$, then with $F_2$, and finally with  $F_3$.  Let $F'_i$ denote the new $f$-surfaces in $N'$  constructed in the manner above. The difference of $F'_1 \cdot F'_2 \cdot {F'_3} - F_1 \cdot F_2 \cdot {F_3}$ is the signed intersection number of the three $f$-surfaces with boundary  in $\nu_{\gamma'}$. This is  $(F_1 \cdot \gamma)(F_2 \cdot \gamma)(F_3 \cdot \gamma)$ times $\tau$, where $\tau$ is  the signed triple-intersection number of three $F$ surfaces with boundary in $\nu_{\gamma'}$ 
which meet the boundary in three parallel curves which are meridians of $\nu_{\gamma}$. However one may easily imagine, in a collar of the boundary, three $f$-surfaces with boundary without any triple intersections  which rearrange the order of these three curves by a single permutation. Since the triple intersection number is skew-symmetric, $\tau$ must be zero
 under the hypothesis  that $f$ is odd.

One may pass from  $S^1 \times S^2$ to the real projective $3$-space  by strong  type-$2$ surgery.  Thus we see that the ring structure on $H^*( \ ,\BZ_2)$  is not preserved by strong $2$-congruence.

\end{proof}

The $f$th Burnside group of a manifold $M$, denoted $B_f(M)$,  is obtained by quotienting the fundamental group of the manifold by the subgroup normally generated by all the $f$th powers of all elements. Dabkowski and  Przytycki's  have considered the $f$-th Burnside group of a double branched cover of a link (this is the $f$th Burnside group of the link)  as a tool in their study of local moves on links. They state in 
\cite[proof of Theorem (1.1)]{DP} that $B_f(M)$ is preserved by $n/f$-surgeries. We note that it is also clear that it is preserved by  $n/sf$ surgeries. Thus one has the following slight generalization of the observation of Dabkowski and  Przytycki.

\begin{prop}\label{burn}
If $M$ and $N$ are weakly $f$-congruent, then $B_f(M)$ and $B_f(N)$ are isomorphic.
\end{prop}

\begin{prop}\label{double} The double branched cover of $S^3$ along a link with  $c$ components  is  strongly $2$-congruent  to the connected sum of  $c-1$ copies of $S^1 \times S^2$.
\end{prop}

This follows from the Montesinos trick \cite{Mo,Li}  that a crossing  change in a link in $S^3$ corresponds to a strong type-2 surgery in the double branched cover of a link. Thus the double  branched cover of $S^3$ along a link with $c$ components is strongly $2$-congruent to the double branched cover of an unlink with $c$ components: the connected sum of   $c-1$ copies of $S^1 \times S^2$.

More generally,  Dabkowski and  Przytycki consider   $\alpha/\beta$-moves between links in $S^3$, where $\alpha$ and $\beta$ are relatively prime integers.  The $\alpha/\beta$-move replaces two parallel strands by a rational tangle specified by $\alpha/\beta$.   Such a move is covered by surgery with numerator $\beta$ and denominator $\alpha$ to the double branched covers of these links  \cite{Mo,DP}. Thus a $fs/n$ move between links is covered by weak type-$f$ surgery with numerator $n$ and denominator $fs$ between their double branched covers.  This is a type-$f$ surgery if and only if $n \equiv \pm 1 \pmod{f}$ and is strong type-$f$ surgery if and only if $n \equiv \pm 1 \pmod{fs}.$ 

\begin{de}
We will  say that a link is rationally $f$-trivial if there is a sequence of $fs/n$ moves ( for varying $s$, and $n$) connecting the link to an unlink.   
\end{de}

This notion of triviality is similar but a somewhat weaker than that considered in \cite{DP}.

\begin{prop}\label{rational} If $L$ is rationally $f$-trivial, then the double branched cover of $S^3$ along a link with  $c$ components  is  weakly $f$-congruent  to the connected sum of  $c-1$ copies of $S^1 \times S^2$.
\end{prop}

\begin{prop}\label{fLink}
Let $N$ (resp. $N'$)  be obtained by Dehn surgery along an ordered link $L$ (resp. $L'$) in $S^3$ described by rational labels  in the manner of Rolfsen. Suppose that $L'$ is obtained from $L$ by a sequence of isotopies and moves  which insert  $f$-full twists  between two  parallel strands of the link. Moreover assume that the labels on $L$  and $L'$ agree modulo $f$, component by component.   Then $N$ and $N'$ are  strongly $f$-congruent.
\end{prop}

\begin{proof} Figure \ref{f2} which is also valid if $n$ is a rational label,  shows how to increment the rational label by $f$. Thus we only need to see how to insert $f$ twists between two strands anywhere one wants  by a strong type -$f$ surgery. But this is by a similar argument.  Perform $-1/f$ along a unknot encircling the two strands and  perform a Rolfsen twist  to twist the strands. The surgery coefficient on the unknot is now $-1/0$, and so it may be erased. The surgery coefficients on the two strands has  gone up by $f$ but this can can be readjusted by a multiple of $f$.
\end{proof}

\begin{figure}[ht!] 
\labellist
    \small
  \hair 2pt
 \pinlabel {$b$} [r] at  0 24.6
  \pinlabel {$a$} [t] at  48.7 5
    \pinlabel {$c$} [l] at  65.7 61.2 
    \pinlabel {$d$} [l] at  97.4 25.1  
      \pinlabel {$a$} [t] at  188.2 5
       \pinlabel {$b$} [r] at  136 24.6
           \pinlabel {$c$} [l] at  203 61.2 
            \pinlabel {$b$} [r] at  263 37
              \pinlabel {$a$} [t] at 318 17
 \endlabellist
\begin{center}
\includegraphics[width=3in]{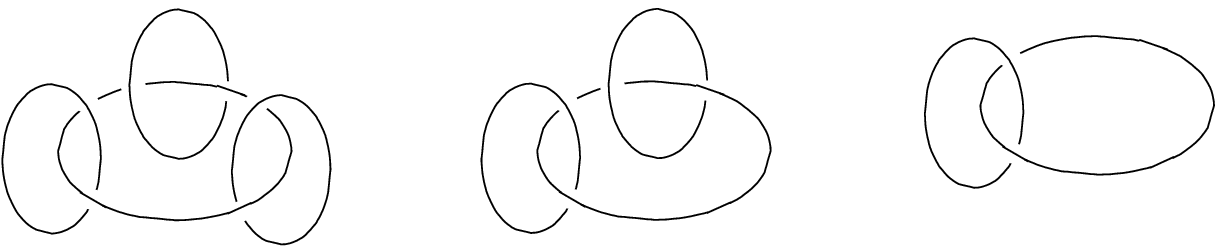}
\caption{The 3-manifolds $H(a,b,c,d)$, $H(a,b,c)$, and $H(a,b)$}
\label{f3}
\end{center}
\end{figure} 

Let $P$ denote the Poincar\'{e} homology  sphere,  and $\Si$ denote  the Brieskorn homology sphere $\Si(2,3,7)$.
\begin{prop}\label{23equiv}
$S^3$, $P$, $-P$, $\Si$ and $-\Si$ are strongly $f$-congruent, 
for $f$ = 2, 3, and 4. $P$ is  strongly $5$-congruent to $-P$. 
 $\Si$ is  strongly $6$-congruent  to $S^3$. 
\end{prop}

\begin{proof} To fix orientations, we take $P$ to be  $-1$ surgery on the left handed trefoil, and $\Si$ to be $-1$ surgery on the right handed trefoil. We use the notations of Figure \ref{f3}. We have that
$P= H(0,-2,3,5)$ and $\Si= H(0,-2,3,7)$. Let $\thickapprox_f$ denote strong $f$-congruence. By Proposition \ref{fLink},
$P= H(0,-2,3,5)\thickapprox_2 H(0,0,1,-1)= H(0,0)= S^3 $, where the  equals  comes from blowing down the $1$ and $-1$. That 
$\Si \thickapprox_2 S^3$ is proved similarly. We have that
$P= H(0,-2,3,5)\thickapprox_3 H(0,1,3,-1)=H(0,3)$, where the  equals  comes from blowing down the $1$ and $-1$. We recognize 
$H(0,3)$ as a genus one homology sphere, i.e. $S^3.$ Similarly, $\Si= H(0,-2,3,7)\thickapprox_3 H(0,1,0,1)= H(2,0)=S^3$. Also
$P= H(0,-2,3,5)\thickapprox_4 H(0,-2,-1,1)= H(0,-2)=S^3$, and $\Si= H(0,-2,3,7)\thickapprox_4 H(0,2,-1,3)=H(1,2,3)= H(1,2)=S^3$. As $S^3=-S^3$, we have obtained the claimed strong $f$-congruences for $f=2, 3, 4$. 

Next $P= H(0,-2,3,5)\thickapprox_5 H(0,-2,-2,0)= L(2,1)\#L(2,1)=  \linebreak L(2,-1)\# L(2,-1)=H(0,2,2,0)\thickapprox_5 H(0,2,-3,-5)=-P$. The identification of $H(0,-2,-2,0)$ holds as one may slide the two components that are framed $-2$ over the fourth component labelled zero, and unlink them from the first component. A zero framed Hopf link yields  $S^3$, and an unknot framed  $-2$ is the lens space $L(2,1)$. Then we make use of the fact that $L(2,1)= L(2,-1)$, as $-1\equiv 1 \pmod{2}$.  Then we slide back over one of the components of the zero framed Hopf  link.

Finally $\Si= H(0,-2,3,7)\thickapprox_6 H(0,-2,3,1)=H(-1,-2,3)= H(-1,4)=U(5) \thickapprox_6 U(-1)=S^3$. Here $U(k)$ denotes $k$ framed surgery along an unknot.
\end{proof}

\begin{prop}\label{Betti}   Let $f$ be a prime. Each $3$-manifold $N$ is strongly $f$-congruent to some 3-manifold $N'$ with $\dim  (H_1(N', \BQ)) = \dim 
(H_1(N', \BZ_f)) = \dim (H_1(N, \BZ_f)).$
\end{prop}

\begin{proof}
Suppose the $N$ is described  as surgery on a framed link $L$ with  $n$ components.
Let $W$ be the result of attaching  2-handles to the $4$-ball according to $L$. Then  $N$ is the boundary of  $W$.  One has that $W$ is simply connected and $H_2(W)= \BZ^n$  with  basis $\{h_i\}$ given by the cores of the  2-handles  capped off in  the  $4$-ball. Similarly 
$H_2(W,N)= \BZ^n$  with basis $\{c_i\}$ given by the  co-cores of the  2-handles. The matrix $\Gamma_L$, associated  to $L$,  with linking numbers on the off-diagonal entries and framings on diagonal entries, is the matrix for the map $H_2(W) \rightarrow H_2(W,N)$ with respect to the bases $\{h_i\}$ and 
$\{c_i\}.$ On the other hand, $\Gamma_L$ is the matrix for the  intersection form $\iota_W$ on $H_2(W)$ with respect to $\{h_i\}$. Consider the commutative diagram:
\[ \begin{CD}  H_2(N)@>{1-1}>>H_2(W) @>>> H_2(W,N)\\
			   	@VVV @V{\pi}VV @VVV     \\
		H_2(N,\BZ_f)@>{1-1}>>H_2(W,\BZ_f) @>{j}>> H_2(W,N,\BZ_f)\
		\end{CD}
\]
The kernel of $j$ has dimension $\dim(H_2(N,\BZ_f))=\dim(H_1(N,\BZ_f))$, which we will denote by $\beta$. We can pick a direct summand $S$ of  
$H_2(W)$ of  dimension $\beta$ which maps to this kernel under $\pi$. Elements $s \in S$ have the property that  for all $x \in H_2(W)$, $\iota_w(x,s) \equiv 0 \pmod{f}.$
Pick a basis  
for $S$, and extend it to a basis $\{\tilde b_i\}$ for $H_2(W)$. Replacing  the first element of this basis by minus itself if necessary we may assume  that the change of basis matrix from $\{b_i\}$ to $\{\tilde b_i\}$ is in $SL(n,\BZ)$ and thus this matrix can be written as a product of elementary matrices: those  with ones on the  diagonal and one $\pm 1$ elsewhere. Changing the basis by an elementary matrix corresponds to a handle slide. Perform a sequence of handle slides  which corresponds to the above product of elementary matrices, and obtain a new framed link $\tilde L$ description of $N$ such that  $\Gamma_{\tilde L}$ has every entry in the  first $\beta$ rows and the  first $\beta$ columns divisible by $p$. Then by Proposition \ref{fLink}, $N$ is  strongly $p$-congruent  to $N'$
where $N'$ is surgery on a framed link $L'$, and $\Gamma_{L'}$ has zero for every entry in the  first $\beta$ rows  and the first $\beta$ columns. By the above exact sequence, but now for  $N'$, we have $\dim(H_2(N', \BQ)) \ge \beta$.  By Theorem \ref{ring},
$\dim(H_2(N', \BZ_f))=\beta$. So by  Poincar\'{e} duality and  the universal coefficients theorem,  $\dim(H_1(N', \BQ)) = \beta =\dim(H_1(N', \BZ_f)) =\dim(H_1(N, \BZ_f)).$ 
\end{proof}

\begin{de} The $f$-cut number of a 3-manifold $N$, denoted $c_f(N)$, is the maximum number of disjoint piecewise linearly embedded good $f$-surfaces that we can place in $N$ with a connected complement.\end{de} 

Recall the  cut number, $c(N)$, is given by the same definition except the surfaces must be  oriented surfaces. One, of course, has
$c_f(N) \ge c(N).$ 

\begin{prop}\label{c_f}   Each $3$-manifold $N$ is weakly $f$-congruent to some 3-manifold $N'$ with $c(N') \ge c_f(N)$. \end{prop}

\begin{proof} Suppose we have $c_f(N)$ disjoint embedded good $f$-surfaces with a connected complement. For each component 
$\gamma$ of the 1-strata of an f-surface $F$, let $\mathcal{T}_\gamma$ be the boundary of a small tubular neighborhood $\nu_\gamma$ of $\gamma$ with $F \cap \nu_\gamma$ consisting of a mapping cylinder for  a covering map $F \cap \mathcal{T}_\gamma \rightarrow \gamma$.   A connected component of $F \cap \mathcal{T}_\gamma$  the represents a multiple of $f$ times the generator for the first homology of $\nu_\gamma$. Each component of $F \cap \mathcal{T}_\gamma$ represents the same homology class of $\mathcal{T}_\gamma$, say $n \mu + fs \lambda$. Thus we may perform a weak type-$f$ surgery along $\gamma$ in such a way that $F \setminus 
\gamma$ may be completed with the addition of one annuli for each component of $F \cap \mathcal{T}_\gamma$ in the surgered manifold. Thus we may perform weak type-$f$ surgery along each component of the 1-strata  in such a way that each the $f$-surface minus a neighborhood of their 1-strata may be completed to an oriented surface. The complement of the resulting the $f$-surfaces remains connected. 
\end{proof}

 \section{Quantum Obstructions to Congruence}
 
 In this paper,  3-manifolds come equipped with a possibly empty $\mathcal{C}$-colored framed link. Here $\mathcal{C}$ is some fixed modular category in the sense of  Bakalov and Kirillov \cite{BK}. Specific manifolds like $S^3$, $S^1 \times S^2$, $P$,  $\Si$,
and $W$ (below) should be assumed to be equipped with the empty framed link unless otherwise stated.  In the definition of  (strong,weak) type-$f$ surgery to a 3-manifold $N$, the surgery curves should be chosen away from $N$'s framed link.  

 Let $\tau_\mathcal{C}(N)\in k_\mathcal{C}$ denote the  Reshetikhin-Turaev quantum invariant  associated to a closed  connected 3-manifold $N$ and the modular category  $\mathcal{C}$ as in \cite[4.1.6]{BK}.  Here $k_\mathcal{C}$ is the ground field of $\mathcal{C}$  algebraically extended, if necessary, so that it contains $D$ and $\zeta$ \cite[3.1.15]{BK}.  The set of isomorphism classes on non-zero simple objects is denoted ${I_\mathcal{C}}$.
 The associated  scalars $\zeta_\mathcal{C}$, and 
 ${\theta_\mathcal{C}}_i$ for  $i \in {I_\mathcal{C}}$ are roots of unity \cite[3.1.19]{BK}. We let $\kappa_\mathcal{C} $ denote 
 ${ \zeta_{\mathcal{C}} }^3$. We say two elements of $k$  whose quotient is a power of 
 $\kappa_\mathcal{C} $ agree up to  phase for $\mathcal{C}$.
 Sometimes we will say simply ``up to phase'', if $\mathcal{C}$ is clear from context.
Let $t_\mathcal{C}$ be the least positive integer $t$ such that 
for some fixed $j \in  \BZ$, ${{\theta_\mathcal{C}}_i}^t= {\kappa_\mathcal{C}}^j$ 
for all for  $i \in {I_\mathcal{C}}$.
Changing the framing by $t_\mathcal{C}$ of any component of a framed link description of a manifold leaves the formula for  $\tau_\mathcal{C}$ unchanged up to phase. 
Thus we have:

\begin{thm}\label{cong}  If $M$ is  strongly $t_\mathcal{C}$-congruent  to $N$, then
$\tau_\mathcal{C}(M)$ and $\tau_\mathcal{C}(N)$ agree up to  phase.
\end{thm}

Let ${\mathcal{V}(n)}$  denote the modular category described by Turaev using the Kauffman bracket skein theory \cite[7.7.1]{T}, with $A$ a primitive $4n$th root of unity, with $n > 3$. 
Then  $\tau_{\mathcal{V}(n)}(M)$ is also known as an $SU(2)$ invariant. It is the same as $\langle  M \rangle_{2n}$ in the notation  of Blanchet, Habegger, Masbaum, and Vogel \cite{BHMV}. 
Moreover $t_\mathcal{C}=4n$, and ${\kappa_{\mathcal{V}(n)}}^2= A^{-6-n(2n+1)}.$ 

\begin{cor}\label{SU}  If $M$ is  strongly $4n$-congruent   to $N$, then $\langle  M \rangle_{2n}$ and $\langle  N \rangle_{2n} $ agree up to  phase for  ${\mathcal{V}(n)}$.
\end{cor}

With the above hypothesis and argument, Lackenby \cite{La} observed the somewhat weaker conclusion: 
$|\langle  M \rangle_{2n}|=  |\langle  N  \rangle_{2n}|.$

{\em  In this paper, $r$ will always denote an odd integer greater than one.} We now consider the modular category ${{\mathcal{V}(r)}^e}$ \cite[7.5]{T}  where the colors of the labels are restricted to be even  integers from $0$  to $r-3$ and $A$ is taken to be a primitive $2rth$ root of unity. \footnote{We make this choice to be consistent with  \cite{BHMV}. Note this is a different choice of $A$  than is made by Tureav. However the quantum invariants of 3-manifolds in both cases are rational functions of $A^4$, and the fourth power of a primitive $4r$th root of unity and of a primitive $2r$th root of unity are both primitive  $r$th roots of unity, and are therefore Galois conjugates. Thus the invariants of closed 3-manifolds so defined differ only by a fixed Galois automorphism.}
Then  $\tau_{{\mathcal{V}(r)}^e}(M)$, also known as the $SO(3)$ invariant, is the same as $\langle  M \rangle_{r}$ in the notation  of \cite{BHMV}. 
Moreover $m_{{\mathcal{V}(r)}^e}=t_{{\mathcal{V}(r)}^e}=r$, and ${\kappa_{{\mathcal{V}(r)}^e}}^2=  A^{-6-r(r+1)/2}.$

We obtain  the following corollary by specializing  Theorem \ref{cong} to the modular category ${{\mathcal{V}(r)}^e}$. 

 \begin{cor}\label{strong}  If $M$ is   strongly $r$-congruent  to $N$, then
$\langle  M  \rangle_r$ and $\langle  N \rangle_r$ agree up to phase for ${\mathcal{V}(r)}^e$.
\end{cor}

\begin{lem}\label{exPBnon}
If there is a strong $r$-congruence  among two of
 $S^3$, $P$, $-P$, $\Si$ and $-\Si$, then one of the following holds:
 \begin{itemize}
\item $r=3$,
\item  $r=5$ and the strong congruence is between $P$ and $-P$, or
\item   $r=7$ the strong congruence is  between $\Si$ and $-\Si$. 
\end{itemize}
\end{lem} 

\begin{proof}  Let $I_r(M)$ denote ${\langle  S^3 \rangle_r}^{-1} \langle  M \rangle_r $

According to Le  \cite{L}, letting  $a$ denote $A^4$ which is a primitive $r$ th root of unity,
\[I_r(P) =  (1 - a)^{-1} \sum_{n=0}^{(r-3)/2} a^n (1 -
a^{n+1}) (1 - a^{n+2}) \cdots (1 - a^{2n+1})
\]
Le has a similar formula for $I_r(\Si)$.

For the case $r >3$ is prime which we denote by $p$, 
we use some techniques of  Chen and Le \cite[\S 6]{CL}. In fact, the argument there together with  our Corollary \ref{strong} show  that the only possible 
strong $p$-congruence (for $p \ge 5$) between $P$ and $-P$ is for $p=5$ and that  the  only possible strong $p$-congruence (for $p \ge 5$) between $\Si$ and $-\Si$ is for $p=7$.  We illustrate 
the method by  showing that $P$ and $\Si$ cannot be strongly $p$-congruent  for $p\ge 5$.  The other stated results are proved in exactly the same way.

Let $h=1-a$. As is well  known, $\BZ[a]/(p)$ is the truncated polynomial ring $\BZ_p[h]/(h^{p-1})$.
If we truncate $a^j$ times the above formula for $I_p(P)$ 
by discarding the terms corresponding to $n>3$ in the index of summation ( which are clearly divisible by $h^4$), and substitute
$a= 1-h$ and then discard terms of order greater than three (in $h$), we obtain a polynomial in $h$ with $\BZ_p$-coefficients which is congruent to the original expression modulo $h^4$. Thus  for primes $p \ge 5$, we have that $a^j I_p(P)$ is congruent  modulo  $h^4$ to
\[  1+(6-j) h+\left(\frac{j^2}{2}-\frac{13 j}{2}+45\right)
   h^2+\left(-\frac{j^3}{6}+\frac{7 j^2}{2}-\frac{145
   j}{3}+464\right) h^3 \]
We used Mathematica \cite{Wo} to work out this and similar expansions. Similarly, using Le's formula for $I_p(\Si)$,  for primes $p >5$, we have that $I_p(\Si)$ is congruent  modulo $h^4$ to
\[ 1+6 h+69 h^2+1064 h^3\]
If $P$ is $p$-congruent  to $\Si$, then, by Corollary \ref{strong},  for some $j$, 
$a^j I_p(P)= \pm I_p(\Si).$ However, noting that both $I_p(P)$, $I_p(\Si)$ are congruent to $1$ modulo $h$, we my discard the $\pm$.
But this implies that the corresponding coefficients
in the two displayed  polynomials in $h$  above must be congruent modulo $p$.
From the coefficients of $h$, $ j \equiv 0 \pmod{p}$. Comparing the coefficients of $h^2$, we conclude that $69-45=24 \equiv 0 \pmod{p}.$ Thus there are no strong $p$-congruences for $p \ge 5.$

Similarly one sees that there can be no strong $p$-congruences for $p>3$, except the
strong $5$-congruence between $P$ and $-P$, and a possible strong $7$-congruence  between $\Si$ and $-\Si.$
In  some of the cases, one must take into account the coefficients of $h^3$. 

Now consider whether there can be a strong $r$-congruence  between $P$ and $\Si$ where $r$ is composite.  Using the contrapositive of Proposition \ref{obv}, the only possibility is that $r$ has the form $3^a.$  We used Mathematica to see that $I_9(P)$ and $I_9(\Si)$ do not agree up to phase. The higher powers of $3$ are then excluded by  Proposition \ref{obv}. The same procedure then works for all pairs of manifolds except the pair $P$ and $-P$, and the pair $\Si$ and $-\Si.$

Proposition \ref{obv} implies the only possible strong $r$-congruences between $P$ and $-P$ are with $r$ of the form $3^a \cdot 5^b$. We used Mathematica to see that $I_9(P)$ and $I_9(\Si)$ disagree up to phase,  that $I_{15}(P)$ and $I_{15}(\Si)$ disagree up to phase, and   that $I_{25}(P)$ and $I_{25}(\Si)$ disagree up to phase. Then by Proposition \ref{obv}, $r$ must be $3$, or $5$.

The pair of  $\Si$ and $-\Si$ is dealt with similarly. \end{proof}

\begin{thm}\label{exPBnon2}
If $P$ and $S^3$ are strongly $f$-congruent, then $f\in \{2,3,4,6,8,12,16, 24\}$. 
If $\Si$ and $S^3$ are strongly $f$-congruent, then $f\in \{2,3,4,6,8,12,16, 24,32\}$.
If $P$ and $\Si$ are strongly $f$-congruent, then $f\in \{2,3,4,6,8,12,16, 24\}$. 
If $P$ and $-\Si$ are strongly $f$-congruent, then $f\in \{2,3,4,6,8,12,16, 24,48\}$.
If $P$ and $-P$ are strongly $f$-congruent, then $f\in \{2,3,4,5,6,8,10,12,16,20,24, 32, 40\}$.
If $\Si$ and $-\Si$ are strongly $f$-congruent, then $f\in \{2,3,4,6,7,8,12,14,16, 24,28, 32, 56\}$.
 \end{thm}
 
 \begin{proof}
Using recoupling theory, one has:
\[ \langle  P  \rangle_{2n}= \kappa \eta^2  \sum_{i=0}^{n-2}
    \theta_i^2  \Delta_i
     \sum_{j=0}^{\min(i,n-2-i)} \Delta_{2j} (\lambda_{2j}^{i\ i})^3\]
     where $\theta_i= (-A)^{i(i+2)}$, $\eta$ and $\kappa$ are as in \cite{BHMV} with ``$p$'' set to $2n$, and $\Delta_i$ and $\lambda_{i}^{j\  k}$ are as in \cite{KL}. Similarly
     \[ \langle  \Si \rangle_{2n}= \kappa \eta^2  \sum_{i=0}^{n-2}
    \theta_i^{-4}  \Delta_i
     \sum_{j=0}^{\min(i,n-2-i)} \Delta_{2j} (\lambda_{2j}^{i\ i})^{-3}.\]
Using Mathematica, we see that $\langle  P \rangle_{16}$ and $\langle  \Si \rangle_{16}$ disagree up to phase. Similarly $\langle  P \rangle_{24}$ and $\langle  \Si \rangle_{24}$ disagree up to phase. By Corollary \ref{SU}, $P$ and $\Si$ cannot be strongly $32$-congruent or strongly $48$-congruent. By Lemma \ref{exPBnon}, there can be no strong $r$-congruences for odd $r>3$.   Using the contrapositive of Proposition \ref{obv}, we have the result for $P$ and $\Si$.

Using Mathematica, we see that $\langle  P \rangle_{32}$ and $\langle  -P \rangle_{32}$ disagree up to phase. Similarly $\langle  P \rangle_{24}$ and $\langle  -P \rangle_{24}$ and also  $\langle  P \rangle_{40}$ and $\langle  -P \rangle_{40}$ disagree up to phase.  By Corollary \ref{SU}, $P$ and $-P$ cannot be strongly  $64$-congruent,  strongly $48$-congruent  or strongly $80$-congruent. By Lemma \ref{exPBnon}, there can be no strong $r$-congruences for odd $r>5$.  Using the contrapositive of Proposition \ref{obv}, we have the result for $P$ and $-P$.

The proofs for other pairs of manifolds are done similarly.
\end{proof}

 Recall that there is an associated projective $SL(2,\BZ)$ action $\rho_\mathcal{C}$ on the $k_\mathcal{C}$ vector space with basis $\mathcal{I}_\mathcal{C}$ where the projective ambiguity is only up to phase, i.e. powers of $\kappa_\mathcal{C}$.

\begin{de}\label{congP} If the representation $\rho_\mathcal{C}$ factors through $SL(2,\BZ_{m})$,  and  $\rho \begin{bmatrix} n^{-1} &0 \\ 0 & n \end{bmatrix}$, with respect to the basis $\mathcal{I}_\mathcal{C}$, is given by  a signed permutation matrix for every invertible $n \in \BZ_{m}$, we will say $\mathcal{C}$ has the $m$-congruence property.
\end{de}

\begin{thm}\label{color}  Suppose that $\mathcal{C}$ satisfies the 
$m$-congruence property.  If $M$ is obtained  from $N$ by a weak type-${m}$ surgery along $\gamma$, then for some color $c \in \mathcal{I}_\mathcal{C}$,  $\tau_\mathcal{C}(M)$ and  $\tau_\mathcal{C}( N \text { with $ \gamma$ colored $c$})$ agree up to phase and sign.
\end{thm}

\begin{proof} We use the associated TQFT.
Our results follow from the  observation that the regluing map for the torus under weak type-${m}$  surgery
can be factored in $SL(2, \BZ_{m})$ as a product of a power of a Dehn  twist on the meridian and  $\begin{bmatrix} n^{-1} & 0 \\ 0 & n \end{bmatrix} $ (see the proof of Theorem \ref{colorprime} below) and that  the representation applied to such a product is particularly simple: the product of a phased permutation matrix and a diagonal matrix which fixes  the unit object  of $\mathcal{I}_\mathcal{C}$.  Note that this factorization  cannot be done in $SL(2, \BZ)$. \end{proof}

According to Bantay \cite[Theorem 3 and equation for $G_\ell$ p.434] {B}, 
modular categories associated to conformal field theories have the $m$-congruence property for some $m$. See also \cite[6.1.7]{Ga}.  We are unsure what the precise mathematical hypotheses are for this result. 

{\em  In this paper, $p$ will always denote an odd prime.}
In the last section, we  study ${{\mathcal{V}(p)}^e}$  and see that it satisfies the $p$-congruence property. According Freedman-Kruskal and independently Larsen-Wang, the associated projective representation factors through an irreducible component of the metaplectic representation of $SL(2, \BZ_p)$ \cite{FK,LW}.   We give the refined version of this result that we need. Our proof is along the lines of \cite{FK}.  An earlier version of this  paper was written only considering the $SO(3)$ theory at odd primes. After learning of Lackenby's earlier work, we placed the  results in the context of  modular categories to highlight the relationship  to \cite{La}.

In Theorem \ref{colorprime} below,  we   derive a more precise version of Theorem \ref{color}  specialized to ${{\mathcal{V}(p)}^e}$ where we specify the color $c$. { \it Let $d$ denote $(p-1)/2$.}  Recall there are $d$ even colors for this theory. We allow links with odd colors as well now as they are allowed in this theory.
Recall an odd colored component may be traded for an even colored component 
\cite[6.3(iii)]{BHMV1}.
  We  found  Theorem \ref{colorprime} surprising as the effect of a general surgery on  the quantum invariant  is the same as replacing the surgery curve by a specified {\em linear combination} of colored curves where the coefficients have {\em denominators}. 

We use the TQFT $(V_p, Z_p)$ of \cite{BHMV} with  $A=-q^d$ and
$\eta = -i (q-q^{-1})/ \sqrt{p}=<S^3>_p$, modified as in \cite{G2} with $p_1$-structures replaced by integral weights on 3-manifolds and lagrangian subspaces of the first homology of surfaces. We use  $\kappa$ to denote the root of unity denoted by $\kappa^3$ in \cite{BHMV}.  Our choice of $A$ and  $\eta$ in section determines $\kappa$ according to the equation \cite[p.897]{BHMV}. Up to sign $\kappa$ is determined by $\kappa^2= A^{-6-p(p+1)/2}$. 
Note that $A^2= q^{-1}$ and so the quantum integers are given by the familiar formula in terms of $q$: $ [n]=\frac {A^{2n} - A^{-2n}}  {A^{2} - A^{-2}} =\frac {q^{n} - q^{-n}}  {q - q^{-1}}$.   We let $\BO_p$ denote the cyclotomic ring of integers $\BZ[A,\kappa_{{\mathcal{V}(p)}^e}]$. If $p=-1 \mod{4}$, $\BO_p$ is $\BZ$ adjoined a primitive $p$th root of unity.  If $p=1 \mod{4}$, $\BO_p$ is $\BZ$ adjoined a primitive $4p$th root of unity.  According to H. Murakami and Masbaum-Roberts \cite{M, MR} $I_p(M) \in \BO_p$. Moreover if  $\beta_1(M)>0$, $\langle  M \rangle_p  \in \BO_p$ \cite{G2}.   If $x/y$ is a unit from $\BO_p$, we write $x  \sim y.$
We have that $\mathcal{D}  \sim (1-q)^{d-1}$.

 \begin{thm}\label{colorprime} Suppose that $M$ is obtained  from $N$ by a weak type-$p$ surgery along $\gamma$ with numerator $n$. Let $\hat n$ denote the integer in the range $[1,d]$ with $n  \hat{n} \equiv \pm 1 \pmod{p}$, and $\check{n}$ denote $\hat{n} -1$. Then for some $m \in \BZ$, 
\[\langle  M \rangle_p =  \kappa^m   \langle   N \text { with $ \gamma$ colored $\check{n}$} \rangle_p  \]
\end{thm}
 \begin{proof} If $p=3$, the quantum invariant for any closed 3-manifold is always $\pm 1$.  Thus  the result holds trivially when $p=3$. So we may assume that  $p\ne 3.$ $M$ is obtained  from $N$ by removing a tubular neighborhood of $\gamma$ and reglueing  by a map $R$ defined by  the matrix 
  \[ \begin{bmatrix} a& ps \\ b &n \end{bmatrix} \equiv \begin{bmatrix} n^{-1} & 0 \\ 0 & n \end{bmatrix} \begin{bmatrix} 1 & 0 \\ 1 &1 \end{bmatrix} ^{n^{-1} b} \equiv U(n) T^{n^{-1} b}\pmod{p}. \]  
  Using the notation of \cite{G2} for vacuum states and pairings,
 by Theorem \ref{factor} and Lemma \ref{un}, $Z(R) [ \nu_\gamma] = \kappa^m  [\nu_\gamma \text{with $\gamma$ colored $\check{n}$}]$, for some integer $m$.  We have 
 \begin{align*}  \langle  M \rangle_p =&\langle  Z[R][ \nu_\gamma] ,[-N \setminus \Int (\nu_\gamma)] \rangle_ {\mathcal T_\gamma}\\
		&=  \kappa^m  \langle  [\nu_\gamma \text{with $\gamma$ colored $\check{n}$}] ,[-N \setminus \Int (\nu_\gamma)] \rangle_ {\mathcal T_\gamma}\\
=& \kappa^m \langle  N \text{ with $\gamma$ colored $\check{n}$}  \rangle_p . 
 \end{align*}
  \end{proof}

In particular, if the numerator is $\pm1 \pmod{p}$, the color $c$ is zero.
Thus we have the following corollary of Theorem \ref{colorprime}
which overlaps with Corollary \ref{strong}. Note that Corollary \ref{strong} requires strong $r$-congruence for $r$ an odd integer, while Corollary \ref{stronger} requires $p$-congruence for $p$ an odd prime.

\begin{cor}\label{stronger} If $M$ is $p$-congruent  to $N$, then
$\langle  M \rangle_p$ and $ \langle   N  \rangle_p $
agree up to  phase for ${\mathcal{V}(p)}^e$.
\end{cor}

\begin{cor}
If there is a $p$-congruence   among two of
 $S^3$, $P$, $-P$, $\Si$ and $-\Si$, then one of the following holds:
 \begin{itemize}
 \item $p=3$,
 \item  $p=5$ and the congruence is between $P$ and $-P$, or
\item   $p=7$ and the congruence is  between $\Si$ and $-\Si$. 
\end{itemize}
\end{cor} 

\begin{proof}This is the same as the $r$ prime case in  the above proof except that one uses Corollary \ref{stronger} instead of Corollary \ref{strong}.
\end{proof}

\begin{rem} One has $\langle \Si \rangle_7 = a^2 \langle -\Si \rangle_7$. So $\Si$ and $-\Si$ satisfy the necessary condition of Corollary \ref{stronger} to be  $7$-congruent. 
One also has  $\langle P \rangle_5 = a^3 \langle -P \rangle_5$. Moreover 
$\langle P \rangle_5 \ne \langle -P \rangle_5$. Given that $P$ and $-P$
are strongly $5$-congruent,  This shows that  Corollaries \ref{strong} and  \ref{stronger} would not be correct if  phase `` up to phase for etc. '' were removed.
\end{rem}

 \begin{prop}\label{Wid} Let $W$ denote the 3-manifold given by 0-framed surgery along the left handed Whitehead link. $W$ may also be described as the double branched cover of $S^3$ along the $(3,6)$ torus link or the Brieskorn manifold $\Si(2,3,6).$ The cohomology  ring of $W$ is the same as that of $\#^2 S^1 \times S^2$ with any coefficients.   \end{prop}

\begin{proof} 
 The double branched cover of $S^3$ along this link is the Brieskorn manifold $\Si(2,3,6)$ \cite[Lemma (1.1)]{Mi}.  Also by Milnor \cite[Theorem(7.1)]{Mi},  $\Si(2,3,6)$ is a circle bundle over a torus with Euler number $-1.$ Then a framed link description of a circle bundle over a torus with Euler number $-1$ is given by the Borromean rings with two components framed zero and one component framed $-1$ 
\cite[Figure 6.1]{GS}. Blowing down the $-1$, we discover that
$\Si(2,3,6)$ is $W$.

$W$ has the same integral cohomology as that of $\#^2 S^1 \times S^2$. In particular, 
 $H^1(W) = \BZ \oplus \BZ$. Thus the trilinear alternating form on $H^1(W)$ must vanish. The trilinear form determines  the rest of the cohomology ring structure using Poincar\'{e} duality. As there is no torsion, the integral cohomology ring determines the cohomology ring with any coefficients. 
 \end{proof}

The following theorem follows  from Dabkowski and  Przytycki's \cite{DP} study of the Burnside group of the $(3,6)$ torus link, and Propositions \ref{burn} and \ref{Wid}. We will give a different proof using quantum invariants.

\begin{thm}
\label{whitehead} For  $p \ge 5$,   $W$ is not weakly $p$-congruent to $\#^2 S^1 \times S^2$.
\end{thm}

\begin{proof}

Using fusion on link strands which meet the 2-sphere factors as well as those that meet the separating 2-sphere (where the connected sum takes place), we have that  $\langle  \#^2 S^1 \times S^2 \text{ with  colored  link} \rangle_p $ must be a multiple of $\eta^{-1}$   which up to phase and units of $\BO$, is $(1-q)^{d-1}.$ Here we trade colors \cite[Lemma 6.3(c)]{BHMV1}, if necessary, so that  the link has only even colors before we perform the above fusion.  By Theorem \ref{colorprime},
 if $W$ were weakly $p$-congruent  to $ \#^2 S^1 \times S^2$, $\langle  W \rangle_p$ would be divisible by $(1-q)^{d-1}$. To complete the proof, we only need to see that this is not the case.

 \begin{figure}[ht!]
 \labellist
    \small
  \hair 2pt
 \pinlabel {$\eta \sum_{k=0}^{d-1} (-1)^k [k+1] \quad $} [r] at  0 90
  \pinlabel {$\omega$} [l] at  101 179 
   \pinlabel {$k$} [r] at  140 15
 \endlabellist

 \begin{center}
\includegraphics[width=1in]{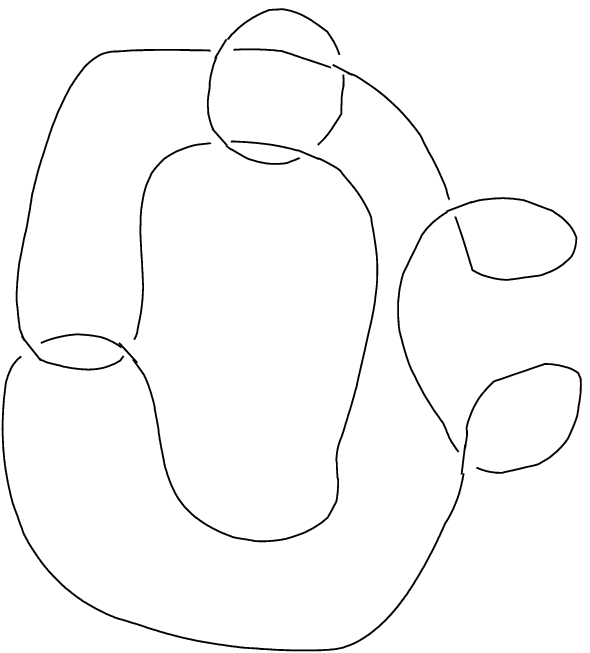}
\caption{Evaluation of $I_p(W)$}
\end{center}
\end{figure}

We wish to calculate first $I_p(W)$ from  Figure 1. Here we give $W$ weight zero.  We apply fusion to the two $k$-colored strands going through the loop colored $\omega$. Only the term with  the strand going through colored zero survives. Moreover the coefficient for this term in the fusion expansion is
$ (-1)^k/ [k+1]$. Also the loop colored $\omega$ with nothing going through it after the fusion contributes 
$\eta^{-1}$. The two positive curls contribute $(-A)^{k(k+2)}$ each, or  $q^{-k(k+2)}$ together. Recall that the evaluation of a Hopf link with both components colored $k$ is $[(k+1)^2]$. Shifting the index of summation, we obtain:
\[ I_p(M) = \sum_{k=0}^{d-1} q^{-k(k+2)} [(k+1)^2] = \frac{q}{q-q^{-1}} \sum_{i=1}^{d} q^{-i^2} (q^{i^2}-q^{-i^2})=
\frac{1}{1-a} \sum_{i=1}^{d}  (1-a^{i^2})\]
where $a=q^{-2}=A^4.$  
As $\eta \sim  (1-a)^{1-d}$, we have  that
 \[\langle  W \rangle_p  = \eta I_p(W) \sim \frac{1}{(1-a)^d} \sum_{i=1}^{d}  (1-a^{i^2}).\]  Thus, summing the same terms twice, adding a zero term,  and using Gauss's quadratic sum,
 \[ 2  (1-a)^d \langle  W \rangle_p \sim \sum_{i=1}^{p-1}  (1-a^{i^2}) = \sum_{i=0}^{p-1}  (1-a^{i^2})\sim  p+\pm i^d \sqrt{p} \]
 As $\sqrt{p} \sim (1-a)^d$, we see that 
 $\langle  W \rangle_p $ is not  a multiple of $1-a$ (or $1-q$).\end{proof}
 
 We now give a slight strengthening of a result of Dabkowski and  Przytycki's \cite[Theorem 2(i)]{DP}. 
 
\begin{cor} The $(3,6)$ torus link is not rationally $p$-trivial
for any prime $p \ge 5.$
\end{cor} 

\begin{proof}  If the $(3,6)$ torus link is  rationally $p$-trivial, then by Propositions \ref{rational} and \ref{Wid}, $W$ would be weakly $p$-congruent to $\#^c S^1 \times S^2,$ for some $c$. As $H_2(W)= \BZ^2, $ we have by Theorem \ref{ring}, that $c$ must be two.
By Theorem \ref{whitehead}, $W$ cannot  be weakly $p$-congruent to $\#^2 S^1 \times S^2$ for $p\ge 5.$
\end{proof}

   \begin{figure}[ht!]
    \labellist
    \small
    \hair 2pt
   \pinlabel $b_1$ [b] at 74 135
       \pinlabel $a$ [l] at 150.1 146.1 
     \pinlabel $c$ [b] at 194.2 163.6 
      \pinlabel $b_2$ [b] at 300 150
     \endlabellist
   \begin{center}
\includegraphics[width=2in]{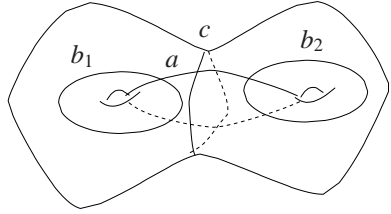}
\caption{Labeling of curves on the boundary of a genus two handlebody} 
\label{handle}
\end{center}
\end{figure} 

Recall that the double of a handlebody of 
genus two $\BH_2$  is the connected sum  of two copies of $S^1 \times S^2$. Let $\BD(h)$ be the result of gluing two copies of $\BH_2$ 
by some element $h$ in the mapping class group of its boundary. If $h$ is in the Torelli group,
then   $\BD(h)$ must have the same cohomology ring as $S^1 \times S^2$ with any coefficients. 
Let \[h= T(b_1)\  
T(b_2)\ T(a)\ T(c) \ T(a)^{-1}\ T(b_2)^{-1}\  T(b_1)^{-1}\]  where $T_x$ denotes a Dehn twist in a neighborhood of a simple closed curve $x$.  As $h$ is conjugate to a Dehn twist around the null-homologous curve $c$, $h$ is in the Torelli  group.  Actually, $\BD(h)=W$.  This identification is an fun exercise in the Kirby calculus making use of the description of the mapping cylinder of a Dehn twist given in Masbaum-Roberts \cite{MR2}.  We can get  variations of Theorem \ref{whitehead} by varying the above word in Dehn twists. For instance, we have:

\begin{thm}\label{ex2}  Let \[h'= T(b_1)\  
T(b_2)\ T(a)^2\ T(c) \ T(a)^{-2}\ T(b_2)^{-1}\  T(b_1)^{-1}.\] The cohomology  ring of 
$ \BD(h')$ is the same as that of $\#^2 S^1 \times S^2$ with any coefficients.  For  $p=5$, $7$, $11$, $13$ and $17$,  $\BD(\text{f})$ is not weakly $p$-congruent  to $\#^2 S^1 \times S^2$. \end{thm}

\begin{proof}  Using  A'Campo's package TQFT  \cite{A'C} which is used with the computer program Pari \cite{PARI}, we have calculated that $\langle  \BD(h') \rangle_p  \in \BO_p$ is not divisible by $1-q$ for the listed $p$.    The rest of the proof is the 
same as the proof of Theorem \ref{whitehead}.\end{proof}

Some other equivalence relations on 3-manifolds which are generated by  surgeries with specified properties were studied by Cochran, Gerges, and Orr \cite{CGO}. The most prominent  of these equivalence relations is called integral homology surgery equivalence. These equivalence relations are different than those
considered here. In particular, $L(f,1)$ is strongly $f$-congruent  to $S^1 \times S^2$
but not  integral homology surgery equivalent. On the other hand, by \cite[Corollary 3.6]{CGO}, $W$ is integral homology surgery equivalent to $\#^2 S^1 \times S^2$  but not weakly $p$-congruent
 to $\#^2 S^1 \times S^2$ for $p \ge 5$ by  {Theorem }\ref{whitehead} .

\section{Quantum invariants}

We now give some results on the integrality and divisibility of quantum invariants.

\begin{thm}\label{cm} If $N$ is a closed   connected 3-manifold with $H_1(N, \BZ_p)$ non-zero,  then  $ \langle  N \rangle_p  \in  \BO_p.$ 
 \end{thm}
 
 \begin{proof}  By Proposition \ref{Betti} $N$ is  strongly $p$-congruent  to a manifold $M$ with positive first Betti number.  By  Corollary   \ref{strong} $\langle  N \rangle_p $ is, up to phase,  the quantum invariant 
of $\langle  M \rangle_p $.  But by \cite[2.12]{G2}  
$\langle  M \rangle_p $ must lie in  $\BO_p.$ 
\end{proof}

 For manifolds without colored links, this is a special case of a result of Cochran and Melvin \cite[Theorem 4.3]{CM} .
This result now holds in the context of 3-manifolds with a  colored link. 
 We also obtain the following  strengthening of  \cite[Theorem 15.1]{GM}

\begin{thm}\label{cpb} If $N$ is a closed   connected 3-manifold with $c_p(N)>0$,  then  \[ \langle  N \rangle_p  \in (1-A^2)^{\frac {(p-3)(c_p(N)-1)} 2 } \BO_p .\] 
 \end{thm}

 \begin{proof}
By Theorem \ref{c_f}, $N$ is weakly $f$-congruent to $N'$ with $c_p(N) \le c(N').$
By repeated use of Theorem \ref{colorprime}, $ \langle  N\rangle_p$ up to phase is given by
 $ \langle  N'  \text {with some colored link} \rangle_p $ which, by \cite[Theorem 15.1]{GM}, 
 is in  $(1-A^2)^{\frac {(p-3)(c(N')-1)} 2 } \BO_p.$
\end{proof}

 \section{$SL(2)$ representations}\label{slrep}
 
 As is well-known,  $SL(2, \BZ)$ is generated  by 
 \[S=\begin{bmatrix} 0 & 1 \\ -1 &0 \end{bmatrix} \quad \text{and} \quad T=\begin{bmatrix} 1 & 0 \\ 1 &1 \end{bmatrix} .\]   They also generate $SL(2, \BZ_p)$ \cite[p.209]{L}.
 
\subsection{Metaplectic representation of  $SL(2, \BZ_p)$}
 Let $d=( p-1)/2$ and $q= e^{\frac{2 \pi i} p}.$  
  We now wish to recall a description of the metaplectic representation of $SL(2, \BZ_p)$ which is given in Neuhauser
 \cite{N}, though the observations about working  over $\BZ[q,1/p]$ are ours. Consider  the $\BC$-vector space $\BC^{\BZ_p}$, the set of complex valued  functions on $\BZ_p.$  It has a basis consisting of the point-characteristic functions $\{\delta_x | x \in  \BZ_p \}$ where $\delta_x(y) =\delta_x^y$ where $x,y \in \BZ_p.$ Using \cite[4.1,4.3,5.6]{N} a true (rather than projective) representation $W$ of   $SL(2, \BZ_p)$ acting  on  $\BC^{\BZ_p}$  can be defined by
 \[ W(S) f(x)= \frac {(-i)^d}{\sqrt{p}} \sum_{y \in \BZ_p} q^{xy} f(y) \quad \text{and} \quad W(T) f(x)= q^{d x^2} f(x)\ . \]
 
 One can see that the prefactor $\frac {(-i)^d}{\sqrt{p}} \in \BZ[q,1/p]$ using  Gauss's quadratic  sum. Thus we can and will view $W$ as a representation on  $\BZ[q,1/p]^{\BZ_p}$. 
For $n \in {\BZ_p}^*$, let $U(n) = \begin{bmatrix} n^{-1} & 0 \\ 0 & n \end{bmatrix} \in SL(2, \BZ_p)$.  By \cite[4.1, 4.3, 5.6]{N}, we have that 
$W (U(n)) f(x) = ( \frac n p )  f(nx).$ Here $( \frac n p )$ denotes the Legendre symbol, and so is $\pm 1.$

The restriction of $W$ to the space of odd functions
 \[\BZ[q,1/p]^{\BZ_p}_\text{odd}=\{f \in \BZ[q,1/p]^{\BZ_p} | f(x)= -f(-x)\}\] is an invariant irreducible summand \cite[4.2]{N} which we denote 
$W_\text{odd}$. 
Define 
$S \subset {\BZ_p}$ by  $S= \{1,2,3, \cdots ,d\}.$ For $x \in S$, define $\delta'_x= \delta_x- \delta_{-x}.$
Then $\{ \delta'_x| x \in S \}$ is a basis for $\BZ[q,1/p]^{\BZ_p}_\text{odd}$ and 
\[\ W_\text{odd}(S) \delta'_x= \frac {(-i)^d}{\sqrt{p}} \sum_{y \in S} (q^{xy}-q^{-xy}) \delta'_y \quad \text{and} \quad W_\text{odd}(T)  \delta'_x  = q^{d x^2} \delta'_x . \]
\begin{lem}\label{un} We have that $W_\text{odd}(U(n))$ sends  $ \delta'_x$ to $ \pm  \delta'_{\pm n^{-1}x}$, where we choose the plus or minus in $\pm n^{-1}x$ so that  $\pm n^{-1}x \in S$. \end{lem}

 \subsection{Projective representation of  $SL(2, \BZ)$ arising from  TQFT}

We use the TQFT $(V_p, Z_p)$ of \cite{BHMV} with  $A=-q^d$ and
$\eta = -i (q-q^{-1})/ \sqrt{p}$, modified as in \cite{G2} with $p_1$-structure replaced by integral weights on 3-manifolds and lagrangian subspaces of the first homology of surfaces.

  We wish to study the projective representation of the mapping class group of the torus ${\mathcal T}$ given by the TQFT $(V_p, Z_p)$. 
  This fails to be an actual representation only by phase factors (powers of $\kappa$).   We think of ${\mathcal T}$ as the boundary of a solid torus 
  ${\mathcal H}$ and 
  pick an ordered basis for the first homology:$ [\lambda]$, $[\mu]$, where $\lambda$ is a longitude and 
  $\mu$ is a meridian. The induced map on homology then defines an isomorphism from the mapping class group of ${\mathcal T}$ to $SL(2, \BZ)$. The map given by  $T$  extends to a full positive twist of the handlebody ${\mathcal H}$. We let $S$ denote the map which is given by $\begin{bmatrix} 0& 1 \\ -1 &0 \end{bmatrix}$.
 This map does not extend over ${\mathcal H}$, but if two copies of ${\mathcal H}$ are glued together using this map (reversing the orientation on the second copy of ${\mathcal H}$), we obtain the 3-sphere with the cores of these handlebodies forming a 0-framed Hopf link.

The module, $V_p({\mathcal T})$,  is  free  with basis $\{e_i | 0 \le i \le d-1 \}$. Here  $e_i$ is  the closure of  the $i$th Jones-Wenzl idempotent in the skein of ${\mathcal H}.$ The basis $b_j= (-1)^{j-1} e_{j-1}$
 for $1 \le j \le d$ is more convenient for us. Using skein theory, one sees that 
 \[  Z(S) b_i = {\eta} \sum_{j=1}^d [ij] b_j \quad \text{and} \quad Z(T) b_i=  (-A)^{i^2-1} b_i\] 
 See,  \cite[p. 2487]{G4} for instance, where the analog is done in the case $p$ is even. Thus
   \[    (-i)^{d-1} Z(S) b_i =  \frac {(-i)^d}{\sqrt{p}} \sum_{j=1}^d (q^{ij}-q^{-ij} )b_j \quad \text{and} \quad  q^{d} Z(T) b_i=  q^{d i^2} b_i\] 
 The factors $(-i)^{d-1}$ and $q^{d}$ are  powers of $\kappa$, except if $p=3$, when $\kappa=-1$.  Modified by  these factors $Z(S)$ and $Z(T)$ are identical to
 $W_{\text{odd}}(S)$ and $W_{\text{odd}}(T)$ under the isomorphism which send $b_i$ to $\delta'_i$. This shows  that the projective representation $Z$ of $SL(2, \BZ)$  can be corrected to an honest representation $Z'$  by rescaling using powers of  
 $\kappa$. 

\begin{thm}\label{factor} Suppose $p \ne 3$. The representation $Z'$ of $SL(2, \BZ)$ factors  through a representation equivalent to the representation $W_\text{odd}$ of  $SL(2, \BZ_p)$  on $\BZ[q, 1/p]^{\BZ_p}_\text{odd}$  under the isomorphism which send  $\delta'_i$ to $b_i$.  Thus $Z$ agrees with $W_\text{odd}$ via this isomorphism, up to powers of $\kappa$, i.e. up to phase.\end{thm}

This refines the result of Freedman-Krushkal and Larsen-Wang  that the projective TQFT representation and the odd part of the projective metaplectic representation are equivalent in $PGL(d, \BC).$

 \begin{rem} The  subgroup $\mathbb{L} $ of $SL(2, \BZ_p)$ consisting of lower triangular matrices is generated by $T$ and the $U(n).$ Thus  $W_\text{odd}$ on this subgroup is represented by ``phased'' permutation matrices. The inverse image of $\mathbb{L}$ in  $SL(2, \BZ)$ consists of the possible glueing matrices for weak type-$p$ surgery.
   \end{rem}

\subsection{Integrality of the metaplectic representation}

In \cite{G2} it is shown that the representation $Z$ of a central extension of $SL(2,Z)$  on $V( \mathcal T)$ lifts to a representation on a module $S( \mathcal T)$ over a cyclotomic ring  of integers. In Gilmer-Masbaum-van Wamelen \cite{GMW}, explicit bases are given in terms of the $e_i$ basis of $V( \mathcal T)$. These results suggest the following.

Let $S$ be the  $\BZ[q]$ submodule of $\BZ[q,1/p]^{\BZ_p}$ generated by  the finite  set
\[ \{ W(g)  \delta_x | g \in SL(2, \BZ_p), x \in \BZ_p \}.\]

\begin{prop}
$S$ is a free finitely generated  $\BZ[q]$ lattice in 
$\BZ[q,1/p]^{\BZ_p}$ of rank $p$ preserved by $SL(2, \BZ_p)$ 
\end{prop}

\begin{proof} We   have that $\BZ[q]$ is a Dedekind domain, $S$ is torsion-free finitely generated  $\BZ[q]$ -module, so   $S$ is  projective \cite{J}. This module becomes free when localized by inverting  $p$. By \cite[lemma 6.2]{G2}, $S$ is already free.
\end{proof}

It would be interesting to find an explicit basis for $S$.

The corresponding results  hold for  $S_\text{odd}$ and $S_\text{even},$ which are defined similarly. In fact one has, by the same proof, the following proposition.

\begin{prop} Let $G$ be a finite group acting on a free finitely generated $\BZ[q,1/p]$-module with basis $\{ \bb_1, \bb_2, \cdots \bb_n\}$. Then  the $\BZ[q]$ submodule generated by 
\[ \{ g \bb_i | g \in G, 1 \le  i \le n \}\]
is a free finitely generated  $\BZ[q]$ lattice of rank  $n$ preserved by 
$G$.
\end{prop}

%
%
%
%

\end{document}